# A note on the $U, V$ method of estimation[*]

## Arthur Cohen[1] and Harold Sackrowitz[1]

*Rutgers University*

**Abstract:** The $U, V$ method of estimation provides unbiased estimators or predictors of random quantities. The method was introduced by Robbins [3] and subsequently studied in a series of papers by Robbins and Zhang. (See Zhang [5].) Practical applications of the method are featured in these papers. We demonstrate that for one $U$ function (one for which there is an important application) the $V$ estimator is inadmissible for a wide class of loss functions. For another important $U$ function the $V$ estimator is admissible for the squared error loss function.

## 1. Introduction

The $U, V$ method of estimation was introduced by Robbins [3]. The method applies to estimating random quantities in an unbiased way, where unbiasedness is defined as follows: The expected value of the estimator equals the expected value of the random quantity to be estimated. More specifically, suppose $X_j$, $j = 1, \ldots, n$, are random variables whose density (or mass) function is denoted by $f_{X_i}(x_i|\theta_i)$. In this paper we consider estimands of the form

$$(1.1) \qquad S(\boldsymbol{X}, \boldsymbol{\theta}) = \sum_{j=1}^{n} U^*(X_j, \theta_j),$$

where $\boldsymbol{X} = (X_1, \ldots, X_n)'$ and $\boldsymbol{\theta} = (\theta_1, \ldots, \theta_n)'$. An estimator, $V(\boldsymbol{X})$ is an unbiased estimator of $S$ if

$$(1.2) \qquad E_{\boldsymbol{\theta}} V(\boldsymbol{X}) = E_{\boldsymbol{\theta}}(S(\boldsymbol{X}, \boldsymbol{\theta})).$$

Of particular interest in applications are estimands of the form $U^*(X_j, \theta_j) = U(X_j)\theta_j$, where $U(\cdot)$ is an indicator function. Robbins [3] offers a number of examples of unbiased estimators using the $U, V$ method. Zhang [5] studies the $U, V$ method for estimating $S$ and provides conditions under which the "$U, V$" estimators are asymptotically efficient. Zhang [5] then presents a Poisson example that deals with a practical problem involving motor vehicle accidents.

In this note we demonstrate that for many practical applications the $U, V$ estimators are inadmissible for many sensible loss functions. In particular, for the Poisson example given in Zhang [5], for the $U$ function given, the $V$ estimator is inadmissible for any reasonable loss function, since the estimator is positive for some $\boldsymbol{X}$ when $S = 0$ no matter which $\theta$ is true.

Previously, Sackrowitz and Samuel-Cahn [4] showed that the $U, V$ estimator of the selected mean of two independent negative exponential distributions is inadmissible for squared error loss.

[*]Research supported by NSF Grant DMS-0457248 and NSA Grant H98230-06-1-007.
[1]Department of Statistics and Biostatistics, Hill Center, Busch Campus, Pisctaway NJ 08854-8019, USA, e-mail: artcohen@rci.rutgers.edu; sackrowi@rci.rutgers.edu

In the next section we examine examples in which $S$ functions based on simple $U$ functions are estimated by inadmissible $V$ functions. For other simple $U$ functions the resulting $V$ estimators are admissible for squared error loss. These later results will be presented in Section 3.

## 2. Inadmissibility results

Let $X_j$, $j = 1, \ldots, n$, be independent random variables with density $f_{X_i}(x_i|\theta_i)$. Let $U^*(X_j, \theta_j) = U(X_j)\theta_j$, where, for some fixed $A \geq 0$,

$$(2.1) \qquad U(X_j) = \begin{cases} 1, & \text{if } X_j \leq A, \\ 0, & \text{if } X_j > A. \end{cases}$$

Consider the following four distributions for $X_j$.

(2.2)        Poisson    $f_X(x|\theta) = e^{-\theta}\theta^x/x!$    $(\theta > 0, x = 0, 1, \ldots)$,

(2.3)        Geometric    $f_X(x|\theta) = (1-\theta)\theta^x$    $(0 < \theta < 1, x = 0, 1, \ldots)$,

(2.4)        Exponential    $f_X(x|\theta) = (1/\theta)e^{-x/\theta}$    $(\theta > 0, x > 0)$,

(2.5)        Uniform Scale    $f_X(x|\theta) = 1/\theta$    $(0 < x < \theta, \theta > 0)$.

Let $W(t)$, $t \geq 0$ be a function with the property that $W(0) = 0$ and $W(t) > 0$ for $t > 0$. Consider loss functions

$$(2.6) \qquad W(a, S) = W(a - S),$$

for action $a$.

For the distributions in (2.2), (2.3), (2.4), (2.5), Robbins [3] finds unique unbiased estimators $V(X_j)$ for $U(X_j)\theta_j$.

**Theorem 2.1.** *Let $X_j$, $j = 1, \ldots, n$, be independent random variables whose distribution is (2.2) or (2.3) or (2.4) or (2.5). Consider the loss function given in (2.6). Let $U(X_j)$ be as in (2.1). Then the unbiased estimator $V(\mathbf{X}) = \sum_{j=1}^n V(X_j)$, where $V(X_j)$ is the unbiased estimator of $U(X_j)\theta_j$, is inadmissible for $S$ given in (1.1).*

*Proof.* The idea of the proof is easily seen if $n = 1$. However for $n > 1$ it is instructive to see how much improvement can be made. The proof for $n = 1$ goes as follows: Let $X_1$ be $X$ and $\theta_1$ be $\theta$. The $V(X)$ estimators for the four cases are given in Robbins [3]. For the Poisson case $V(X) = U(X-1)X$ $(V(0) = 0)$. Now let $[A]$ denote the largest integer in $A$ less that $A$. Then $V([A]+1) = [A]+1$, whereas $S = U([A]+1)\theta = 0$.

If

$$V^*(X) = \begin{cases} V(X), & \text{all } X \text{ except } X = [A]+1, \\ 0, & X = [A]+1, \end{cases}$$

then clearly $V^*(X)$ is better than $V(X)$ since $W(V^*([A]+1) - S) = 0$ for $V^*$ and $W(([A]+1) - S) > 0$ for $V$. For the case of arbitrary $n$, $S = 0$ whenever all $X_j \geq ([A]+1)$ whereas $V(X) \neq 0$ whenever at least one $X_j = ([A]+1)$. If all $X_j = ([A]+1)$, then $V = n([A]+1)$. Clearly if $V^* = 0$ at such $\mathbf{X}$, $V^*$ is better than $V$.

For the geometric distribution when $n = 1$, $V(X) = \sum_{i=0}^{X-1} U(i)$ $(V(0) = 0)$. Note $S = 0$ for $X \geq [A]+1$ but $V = [A]+1$ for all such $X$. Again if $V^* = V$ for $X \leq [A]$ and $V^* = 0$ for $X \geq [A]+1$, $V^*$ is better than $V$. The case of arbitrary $n$ is even more dramatic than is the Poisson case with $S = 0$ if all $X_j \geq [A] + 1$ whereas $V \neq 0$ on such points.

TABLE 1
*Improvement in risk for squared error loss function*

| $A$ \ $n$ | 1 | 2 | 3 | 4 | 5 | 6 | 7 | 8 | 9 | 10 |
|---|---|---|---|---|---|---|---|---|---|---|
| 1 | 1.083 | 1.872 | 2.190 | 2.148 | 1.902 | 1.575 | 1.243 | 0.947 | 0.701 | 0.508 |
| 3 | 3.126 | 4.763 | 5.086 | 4.626 | 3.831 | 2.982 | 2.220 | 1.599 | 1.122 | 0.771 |
| 5 | 5.782 | 8.268 | 8.419 | 7.364 | 5.894 | 4.447 | 3.216 | 2.253 | 1.539 | 1.031 |
| 7 | 8.934 | 12.268 | 12.113 | 10.328 | 8.083 | 5.976 | 4.242 | 2.919 | 1.961 | 1.292 |
| 9 | 12.511 | 16.694 | 16.120 | 13.490 | 10.388 | 7.568 | 5.299 | 3.600 | 2.389 | 1.556 |

For the exponential distribution when $n = 1$, $V(X) = \int_0^X U(t)dt = X$ if $X \leq A$, and $V(X) = A$ if $X > A$. For arbitrary $n$, $S = 0$ whenever all $X_j > A$, whereas $V(\mathbf{X}) \neq 0$ on such points.

For the scale parameter of a uniform distribution with $n = 1$, $V(X) = XU(X) + \int_0^X U(t)dt$ which becomes $2X$ if $X \leq A$ and $A$ if $X > A$. Hence as in the previous case, for arbitrary $n$, $S = 0$ whenever all $X_j > A$ whereas $V(\mathbf{X}) \neq 0$ on such points. This completes the proof of the theorem. $\square$

**Remark 2.1.** *Theorem 2.1 applies to the Poisson example in Zhang [5].*

**Remark 2.2.** *If the loss function in (2.6) is squared error then the amount of improvement in risk of $V^*$ over $V$ depends on $n$, $A$, and $\boldsymbol{\theta}$. It can be easily calculated. For the case where all the components of $\boldsymbol{\theta}$ are equal and each $\theta_i$, $i = 1, \ldots, n$ is set equal to $[A] + 1$ the amount of improvement is equal to*

$$\frac{\sum_{i=1}^n \left(i([A]+1)\right)^2 C_i^n e^{-([A]+1)}([A]+1)^{[A]+1}}{([A]+1)!}$$
(2.7)
$$\cdot \left(\frac{1 - \sum_{y=0}^{[A]+1} e^{-([A]+1)}([A]+1)^y}{y!}\right)$$

Table 1 offers the amount of improvement for $n = 1(1)10$ and for values of $A = 1, 3, 5, 7, 9$. We observe as $n$ gets large the amount of improvement becomes smaller. Also for small $n$ as $A$ gets large, improvement gets large. Such observations are consistent with the asymptotic efficiency of the $U, V$ estimator as $n \to \infty$ and with Sterling's formula.

**Remark 2.3.** *Theorem 2.1 also holds for predicting*

$$S^* = \sum_{j=1}^n Y_j U(X_j),$$

*where $Y_j$ has the same distribution of $X_j$ but is unobserved.*

## 3. Admissibility results

In this section we consider the case

(3.1) $$U(X_j) = \begin{cases} 0, & \text{if } X_j \leq A, \\ 1, & \text{if } X_j > A, \end{cases} \quad A \geq 0; \; j = 1, \ldots, n.$$

Also we consider a squared error loss function.

**Theorem 3.1.** *Suppose $X_j$ are independent with Poisson distributions with parameter $\lambda_j$. Then $V(\boldsymbol{X})$ is an admissible estimator of $S(\boldsymbol{X}, \boldsymbol{\lambda})$ for squared error loss.*

*Proof.* Let $n = 1$ and recall $V(X_1) = U(X_1 - 1)X_1$, $V(0) = 0$. Then

$$V(X) = \begin{cases} 0, & \text{for } X_1 = 0, 1, \ldots, [A] + 1, \\ X_1, & \text{for } X_1 > [A] + 1, \end{cases}$$

while

$$U^*(X_1, \lambda_1) = U(X_1)\lambda_1 = \begin{cases} 0, & X_1 \leq [A], \\ \lambda_1, & X_1 \geq [A] + 1 \end{cases}$$

Since $U^*(X_1, \lambda_1) = 0$ for $X_1 \leq [A]$, any admissible estimator of $U^*(X_1, \lambda_1)$ must estimate 0 for $X_1 \leq [A]$ as $V(X_1)$ does. □

At this point we can restrict the class of estimators to all those which estimate by the value 0 for all $X_1 \leq [A]$. For $[X_1] \geq [A] + 1$, $U^*(X_1, \lambda_1) = \lambda_1$ and we have a traditional problem of estimating a parameter $\lambda_1$. Now we can refer to the proof of Lemma 5.2 of Brown and Farrell [1] to conclude that any estimator that can beat $V(X)$ would have to estimate 0 at $X_1 = [A] + 1$. Furthermore for the conditional problem given $X_1 > [A] + 1$, it follows by results in Johnstone [2] that $X_1$ is an admissible estimator of $\lambda_1$.

For arbitrary $n$ the proof is more detailed. We give the details for $n = 2$. The extension for arbitrary $n$ will follow the steps for $n = 2$ and employ induction. For $n = 2$, suppose $V(X_1) + V(X_2)$ is inadmissible. Then there exists $\delta^*(X_1, X_2)$ such that

$$\sum_{x_1=0}^{\infty} \sum_{x_2=0}^{\infty} \left(V(x_1) + V(x_2) - U(x_1)\lambda_1 - U(x_2)\lambda_2\right)^2 \lambda_1^{x_1} \lambda_2^{x_2} e^{-\lambda_1 - \lambda_2} / x_1! x_2!$$

$$(3.2) \geq \sum_{x_1=0}^{\infty} \sum_{x_2=0}^{\infty} \left(\delta^*(x_1, x_2) - U(x_1)\lambda_1 - U(x_2)\lambda_2\right)^2 \lambda_1^{x_1} \lambda_2^{x_2} e^{-\lambda_1 - \lambda_2} / x_1! x_2!$$

for all $\lambda_1 > 0$, $\lambda_2 > 0$, with strict inequality for some $\lambda_1$ and $\lambda_2$.

Now let $\lambda_2 \to 0$. Then by continuity of the risk function, (3.2) leads to

$$(3.3) \quad E\left\{\left(V(X_1) - U(X_1)\lambda_1\right)^2\right\} \geq E\left\{\left(\delta^*(X_1, 0) - U(X_1)\lambda_1\right)^2\right\}.$$

Since $V(X_1)$ is admissible for $U(X_1)\lambda_1$, the case $n = 1$, (3.3) implies that $V(X_1) = \delta^*(X_1, 0)$. At this point we do as in Brown and Farrell [1] by dividing both sides of (3.2) by $\lambda_2$. Reconsider (3.2) but now we can let the sum on $x_2$ run from 1 to $\infty$ since $V(X_1) = \delta^*(X_1, 0)$. Again let $\lambda_2 \to 0$ and this leads to $V(X_1) = \delta^*(X_1, 1)$. Repeat the process for $X_2 = 0, 1, \ldots, [A] + 1$. Furthermore by symmetry $V(X_2) = \delta^*(0, X_2) = \cdots = \delta^*([A] + 1, X_2)$. Thus $V(X_1) + V(X_2) = \delta^*(X_1, X_2)$ on all sample points except the set $B = (X_1 \geq [A] + 2, X_2 \geq [A] + 2)$. Here $V(X_1) + V(X_2) = X_1 + X_2$ and $S = \lambda_1 + \lambda_2$. We consider the conditional problem of estimating $\lambda_1 + \lambda_2$ by $X_1 + X_2$ given $\boldsymbol{X} \in B$. Clearly when $\lambda_1 = \lambda_2 = \lambda$ no estimator can match, much less beat the risk of $X_1 + X_2$ for this conditional problem since $X_1 + X_2$ is a sufficient statistic, the loss is squared error, and $X_1 + X_2$ is an admissible estimator of $2\lambda$. Thus $\delta^*(X_1, X_2) = V(X_1) + V(X_2)$ on the entire sample space proving the theorem.